\documentclass[reqno]{amsart}
\usepackage{amssymb,setspace}
\usepackage{ifpdf}
\ifpdf
 \usepackage[hyperindex,pagebackref]{hyperref}
\else
 \expandafter\ifx\csname dvipdfm\endcsname\relax
 \usepackage[hypertex,hyperindex,pagebackref]{hyperref}
 \else
 \usepackage[dvipdfm,hyperindex,pagebackref]{hyperref}
 \fi
\fi
\allowdisplaybreaks[4]
\numberwithin{equation}{section}
\theoremstyle{plain}
\newtheorem{thm}{Theorem}[section]

\theoremstyle{remark}

\begin{document}

\title[Limits for ratios of gamma functions at singularities]
{Limit formulas for ratios of derivatives of the gamma and digamma functions at their singularities}

\author[F. Qi]{Feng Qi}
\address{Department of Mathematics, School of Science, Tianjin Polytechnic University, Tianjin City, 300387, China; School of Mathematics and Informatics, Henan Polytechnic University, Jiaozuo City, Henan Province, 454010, China}
\email{\href{mailto: F. Qi <qifeng618@gmail.com>}{qifeng618@gmail.com}, \href{mailto: F. Qi <qifeng618@hotmail.com>}{qifeng618@hotmail.com}, \href{mailto: F. Qi <qifeng618@qq.com>}{qifeng618@qq.com}}
\urladdr{\url{http://qifeng618.wordpress.com}}



\begin{abstract}
In the paper the author simply presents the limit formulas for ratios of derivatives of the gamma function and the digamma function at their singularities.
\end{abstract}

\thanks{This paper was typeset using \AmS-\LaTeX}

\maketitle

\section{Introduction}

Throughout this paper, we use $\mathbb{N}$ to denote the set of all positive integers.
\par
It is well known that the gamma function $\Gamma(z)$ is single valued and analytic over the entire complex plane, save for the points $z=-n$, with $n\in\{0\}\cup\mathbb{N}$, where it possesses simple poles with residue $\frac{(-1)^n}{n!}$. Its reciprocal $\frac1{\Gamma(z)}$ is an entire function possessing simple zeros at the points $z=-n$, with $n\in\{0\}\cup\mathbb{N}$. See related texts in~\cite[p.~255, 6.1.3]{abram}. This implies that
\begin{equation}\label{Gamma(z)-f=n(z)}
\Gamma(z)=\frac{(-1)^{n}}{n!(z+n)}f_n(z)
\end{equation}
is valid on the neighbourhood
\begin{equation}\label{neibour-D}
D\biggl(-n,\frac14\biggr)=\biggl\{z:|z+n|<\frac14\biggr\}
\end{equation}
of the points $z=-n$ with $n\in\{0\}\cup\mathbb{N}$, where $f(z)$ is analytic on $D\bigl(-n,\frac14\bigr)$ and satisfies
\begin{equation}\label{f-n(z)-lim}
\lim_{z\to-n}f_n(z)=1
\end{equation}
for all $n\in\{0\}\cup\mathbb{N}$.
\par
It is also well known that the polygamma functions are defined by $\psi(z)=\frac{\Gamma'(z)}{\Gamma(z)}$ and $\psi^{(i)}(z)$ for $i\in\mathbb{N}$. Among them, the first five functions $\psi(z)$, $\psi'(z)$, $\psi''(z)$, $\psi^{(3)}(z)$, and $\psi^{(4)}(z)$ are known as the di-, tri-, tetra-, pentra-, and hexa-gamma functions respectively.
The polygamma function $\psi^{(n)}(z)$ for $n\in\{0\}\cup\mathbb{N}$ is single valued and analytic over the entire complex plane, save at the points $z=-m$, with $m\in\{0\}\cup\mathbb{N}$, where it possesses poles of order $n+1$. See related texts in~\cite[p.~260, 6.4.1]{abram}. From~\eqref{Gamma(z)-f=n(z)}, it follows that the expressions
\begin{equation}\label{psi=n=sing}
\psi^{(n)}(z)=\frac{(-1)^{n+1}n!}{(z+m)^{n+1}}+\biggl[\frac{f_m'(z)}{f_m(z)}\biggr]^{(n)}
\end{equation}
for $n\in\{0\}\cup\mathbb{N}$ are valid on $D\bigl(-m,\frac14\bigr)$ which is defined by~\eqref{neibour-D}.
\par
In~\cite{Prabhu-arxiv, Prabhu-Srivastava}, by Euler's reflection formulas for $\Gamma(z)$ and $\psi(z)$,  the limit formulas
\begin{equation}\label{gamma-limit-eq}
\lim_{z\to-k}\frac{\Gamma(nz)}{\Gamma(qz)}=(-1)^{(n-q)k}\frac{q}{n}\cdot\frac{(qk)!}{(nk)!}
\end{equation}
and
\begin{equation}\label{polygamma-limit-eq}
  \lim_{z\to-k}\frac{\psi(nz)}{\psi(qz)}=\frac{q}{n}
\end{equation}
for any non-negative integer $k$ and all positive integers $n$ and $q$ were established.
\par
In~\cite{polygamma-sigularity.tex}, by using the explicit formulas for the $n$-th derivatives of the cotangent functions in~\cite{derivative-tan-cot.tex}, the limit formulas
\begin{equation}\label{deriv-polygamma-singul-lim}
\lim_{z\to-k}\frac{\psi^{(i)}(nz)}{\psi^{(i)}(qz)}= \biggl(\dfrac{q}{n}\biggr)^{i+1}
\end{equation}
for $n,q\in\mathbb{N}$ and $i,k\in\{0\}\cup\mathbb{N}$ were presented. It is clear that the limit~\eqref{deriv-polygamma-singul-lim} for $i=0$ becomes~\eqref{polygamma-limit-eq}.
\par
The aim of this paper is to discover the limit formulas for ratios of derivatives of the gamma functions at their singularities.
\par
Our main result may be stated as the following theorem.

\begin{thm}\label{gamma-der-singularity-thm}
For $n,q\in\mathbb{N}$ and $i,k\in\{0\}\cup\mathbb{N}$, we have
\begin{equation}\label{deriv-gamma-singul-lim}
\lim_{z\to-k}\frac{\Gamma^{(i)}(nz)}{\Gamma^{(i)}(qz)} =(-1)^{(n-q)k}\biggl(\frac{q}{n}\biggr)^{i+1}\frac{(qk)!}{(nk)!}.
\end{equation}
\end{thm}

Finally, we provide a very simple proof of the formula~\eqref{deriv-polygamma-singul-lim}.

\section{Proof of Theorem~\ref{gamma-der-singularity-thm}}

Now we set off to prove Theorem~\ref{gamma-der-singularity-thm}.
\par
When $i=0$, the limit~\eqref{deriv-gamma-singul-lim} becomes~\eqref{gamma-limit-eq}.
\par
Differentiating $i\ge0$ times on both sides of~\eqref{Gamma(z)-f=n(z)} yields
\begin{align*}
\Gamma^{(i)}(z)&=\frac{(-1)^{n}}{n!}\sum_{\ell=0}^i\binom{i}{\ell} \biggl(\frac{1}{z+n}\biggr)^{(\ell)}f_n^{(i-\ell)}(z)\\
&=\frac{(-1)^{n}}{n!}\sum_{\ell=0}^i\binom{i}{\ell} \frac{(-1)^\ell\ell!}{(z+n)^{\ell+1}}f_n^{(i-\ell)}(z).
\end{align*}
Therefore, we have
\begin{align*}
\lim_{z\to-k}\frac{\Gamma^{(i)}(nz)}{\Gamma^{(i)}(qz)}&= \lim_{z\to-k}\frac{\frac{(-1)^{nk}}{(nk)!}\sum_{\ell=0}^i\binom{i}{\ell} \frac{(-1)^\ell\ell!}{(nz+nk)^{\ell+1}}f_{nk}^{(i-\ell)}(nz)} {\frac{(-1)^{qk}}{(qk)!}\sum_{\ell=0}^i\binom{i}{\ell} \frac{(-1)^\ell\ell!}{(qz+qk)^{\ell+1}}f_{qk}^{(i-\ell)}(qz)}\\
&=(-1)^{(n-q)k}\biggl(\frac{q}{n}\biggr)^{i+1}\frac{(qk)!}{(nk)!}.
\end{align*}
The proof of Theorem~\ref{gamma-der-singularity-thm} is completed.

\section{A simple proof of the formula~\eqref{deriv-polygamma-singul-lim}}

In virtue of~\eqref{psi=n=sing}, we have
\begin{equation*}
\lim_{z\to-k}\frac{\psi^{(i)}(nz)}{\psi^{(i)}(qz)}=\lim_{z\to-k} \frac{\frac{(-1)^{i+1}i!}{(nz+nk)^{i+1}}+\bigl[\frac{f_{nk}'(nz)}{f_{nk}(nz)}\bigr]^{(i)}} {\frac{(-1)^{i+1}i!}{(qz+qk)^{i+1}}+\bigl[\frac{f_{qk}'(qz)}{f_{qk}(qz)}\bigr]^{(i)}}
= \biggl(\dfrac{q}{n}\biggr)^{i+1}.
\end{equation*}
The proof of the formula~\eqref{deriv-polygamma-singul-lim} is completed.

\end{document}